\documentclass[a4paper,12pt]{article}
\usepackage[T2A]{fontenc}
\usepackage[utf8]{inputenc} 
\usepackage{amssymb,amsfonts}
\usepackage{amsfonts,amsmath,amssymb,amscd,latexsym}
\usepackage{srcltx}

\linethickness{0.4pt}

\newcommand{\bi}{\bibitem}
\newcommand{\nb}{\newblock}

\newcommand{\be}[1]{\begin{equation}\label{#1}}
\newcommand{\ee}{\end{equation}}

\newcommand{\la}{\langle\,}
\newcommand{\ra}{\,\rangle}

\newcommand{\prf}{{\bf Proof.}\ }

\newcommand{\bdelta}{\bar\delta}
\newcommand{\distan}{\mathop{\rm dist}}

\newcommand{\hgt}{\mathop{\rm ht}}

\newtheorem{thm}{\quad Theorem}

\newtheorem{prop}{\quad Proposition}
\newtheorem{rk}{\quad Remark}
\newtheorem{conj}{\quad Conjecture}

\title{Evacuation schemes on Cayley graphs and non-amenability of groups}
\author{\vspace{2ex}
Victor Guba\thanks{This work is supported by the Russian Foundation
for Basic Research, project no. 20-01-00465.}\\
Vologda State University,\\
15 Lenin Street,\\
Vologda\\
Russia\\
160600\\
E-mail: gubavs{@}vogu35.ru}
\date{}

\begin{document}

\maketitle

\begin{abstract}

In this paper we introduce a concept of an evacuation scheme on the Cayley graph of an infinite finitely generated group. This is a collection of infinite simple paths bringing all vertices to infinity. We impose a restriction that every edge can be used a uniformly bounded number of times in this scheme. An easy observation shows that existing of such a scheme is equivalent to non-amenability of the group. A special case happens if every edge can be used only once. These scheme are called pure. We obtain a criterion for existing of such a scheme in terms of isoperimetric constant of the graph. We analyze R.\,Thompson's group $F$, for which the amenability property is a famous open problem. We show that pure evacuation schemes do not exist for the set of generators $\{x_0,x_1,\bar{x}_1\}$, where $\bar{x}_1=x_1x_0^{-1}$. However, the question becomes open if edges with labels $x_0^{\pm1}$ can be used twice. Existing of pure evacuation scheme for this version is implied by some natural conjectures.

\end{abstract}
\vspace{5ex}

Here is the plan of the paper. In Section~\ref{intr} we recall some basic information on the subject: Cayley graphs of groups, isoperimetric constants, (non)amenability, flows on graphs, rooted binary trees (forests), etc. In Section~\ref{ev} we introduce the main concept: an evacuation scheme on a Cayley graph of a finitely-generated group. Roughly speaking, this means that to any vertex we assign an infinite path to bring this vertex to infinity. We impose a natural restriction that evacuation paths can use the same edge only finitely many times, uniformly bounded. It turns out that existence of such a scheme on the Cayley graph of a group is equivalent to its non-amenability.

The case of our special interest happens if edges can be used only once in the family of evacuation paths. This is called a pure evacuation scheme. We find a criterion for existing such a scheme on the Cayley graph: it exists if and only if the Cheeger isoperimetric constant is at least $1$. This is equivalent to the following: for any finite subgraph $Y$, the number of external edges (connecting a vertex in $Y$ to a vertex outside $Y$) is at least $|Y|$. This is a generalization of our previous result for groups with two generators.

Besides, we show that existing of evacuation schemes on the whole Cayley graph is equivalent to existing of such schemes on all finite subgraphs. In this case each vertex goes to the boundary of $Y$ instead of going to infinity. This is also applied to pure evacuation schemes.

Section~\ref{isop} concerns the problem whether Thompson's group $F$ is amenable. We analyze known estimates for various systems of its generators, including Brown -- Belk construction from~\cite{Be04,BB05}. In Theorem~\ref{isopsym0} we strengthen results from~\cite{Gu21a} showing that the isoperimetric constant for the set of generators $\{x_1,\bar{x}_1=x_1x_0^{-1},x_0\}$ is strictly less than $1$. Thus there are no pure evacuation schemes on the corresponding Cayley graph. However, if we slightly change this generating set adding a second copy of generator $x_0$, then we show that some natural conjectures imply existing of a pure evacuation scheme on the Cayley graph for a multi-set $\{x_1,\bar{x}_1,x_0,x_0\}$. So we have an alternative for the amenability problem of $F$: either the above scheme does exist, which implies non-amenability, or some well-confirmed conjectures fail, which will give a breakthrough towards proving amenability.

\section{Introductory concepts}
\label{intr}

\subsection{Cayley graphs}
\label{cay}

By a {\em graph} we mean a (non-oriented) graph in the sense of Serre~\cite{Se80}. This means that each geometric edge consists of two mutually inverse directed edges. As a formal concept, this is a 5-tuple $\Gamma=\la V,E,^{-1},\iota,\tau\ra$, where $V$, $E$ are disjoint sets (of {\em vertices} and {\em edges}, respectively); $^{-1}\colon E\to E$, $\iota\colon E\to V$, $\tau\colon E\to V$. We assume that $e\ne e^{-1}$, $(e^{-1})^{-1}=e$, $\iota(e^{-1})=\tau(e)$, $\tau(e^{-1})=\iota(e)$ for all $e\in E$. Here $e^{-1}$ is called the {\em inverse} edge of $e$, $\iota(e)$ is the {\em initial} vertex of $e$, $\tau(e)$ is the {\em terminal} vertex of $e$.

A {\em path} in a graph $\Gamma$ is defined in a standard way. This is either a single vertex (called a {\em trivial} path), or a sequence of edges written as $p=e_1\ldots e_n$, where $\tau(e_i)=\iota(e_{i+1})$ for all $1\le i < n$. We also need {\em infinite paths} written as $p=e_1e_2\ldots e_n\ldots$ with the same condition $\tau(e_i)=\iota(e_{i+1})$ for all $i$.

Let $A$ be an {\em alphabet}, that is, a set of symbols. Taking a disjoint bijective copy $A^{-1}$ of the set $A$, we get the set $A^{\pm1}=A\cup A^{-1}$ called a {\em group alphabet}. It has an involution $^{-1}$ without fixed points defined by $(a^{-1})^{-1}\rightleftharpoons a$. Elements of the free monoid $(A^{\pm1})^{\ast}$ are called {\em group words}.

Let $G$ be a group equipped by mapping $A\to G$ such that the image of $A$ generates $G$. We will say that $A$ is a set of (group) generators for $G$. Usually $A$ is identified with its image provided the above map is injective. For our needs it is convenient to include the case when different symbols may denote the same group generator. In this case we may say that $A$ defines a multiset of generators, say, $\{a_1=x_0,a_2=x_0,a_3=x_1\}$, where $x_0$ is doubly repeated.

Let $G$ be a group generated by $A$ in the above sense. Its {\em right Cayley graph} $\Gamma_r={\cal C}(G;A)$ is defined as follows. The set of vertices is $G$; for any $a\in A^{\pm1}$ we put a directed edge $e=(g,a)$ with $\iota(e)=g$, $\tau(e)=ga$. Its inverse is $e^{-1}=(ga,a^{-1})$. Each edge has a {\em label} denoted by $\phi(e)$. For the above edge $e=(g,a)$ it is just $a$. The label function can be naturally continued to the set of paths. Namely, for $p=e_1\ldots e_n$ we set $\phi(p)\rightleftharpoons\phi(e_1)\ldots\phi(e_n)$ having a group word over $A$. The label of a trivial path is an empty word.

The concept of a {\em left Cayley graph} $\Gamma_l={\cal C}(G;A)$ is defined in a similar way. Here the set of vertices is also $G$. A directed egde $e=(a,g)$ with label $a$ has $\iota(e)=ag$, $\tau(e)=g$. That is, going along the edge labelled by $a$, means cancelling $a$ on the left. The inverse egde is $e^{-1}=(a^{-1},ag)$. Labels of paths are defined similarly for this case.

\subsection{Automata and isoperimetric constants}
\label{aut}

The cardinality of a finite set $Y$ will be denoted by $|Y|$. 

Let $G$ be an infinite group generated by $A$. For our needs we assume that $A$ is always finite, $|A|=m$.  Let $\Gamma={\cal C}(G,A)$ be the Cayley graph of $G$, right or left. To any finite nonempty subset $Y\subset G$ we assign a subgraph in $\Gamma$ adding all edges connecting vertices of $Y$. So given a set $Y$, we will usually mean the corresponding subgraph. This is a labelled graph that we often call an {\em automaton}. For each $g\in Y$ we have exactly $2m$ directed edges in $\Gamma$ starting at $g$, where $a\in A^{\pm1}$. If the endpoint of such an edge with label $a$ belongs to $Y$, then we say that the vertex $g$ of our automaton $Y$ {\em accepts} $a$. For the case of right Cayley graphs this means $ga\in Y$, for the case of left Cayley graphs this means $a^{-1}g\in Y$. 

A vertex $g\in Y$ is called {\em internal} whenever it accepts all labels $a\in A^{\pm1}$. That is, the degree of $g$ in $Y$ equals $2m$. Otherwise we say that $g$ belongs to the {\em inner boundary} of $Y$ denoted by $\partial A$. 

By $\distan(u,v)$ we denote the distance between two vertices in $\Gamma$, that is, the length of a shortest
path in $\Gamma$ that connects vertices $u$, $v$. For any vertex $v$ and a number $r\ge0$ let $B_r(v)$ denote the ball of radius $r$ around $v$, that is, the set of all vertices in $\Gamma$ at distance $\le r$ from $v$. For any set $Y$ of vertices, by $B_r(Y)$ we denote the $r$-neighbourhood of $Y$, that is, the union of all balls $B_r(v)$, where $v$ runs over $Y$. By
$\partial_o Y$ we denote the {\em outer boundary\/} of $Y$, that is, the set $B_1(Y)\setminus Y$.

An edge $e$ is called {\em internal} whenever it connect two vertices of $Y$. If an edge $e$ connects a vertex of $Y$ with a vertex outside $Y$, then we call it {\em external}. That is, $e$ connects a vertex in $\partial Y$ with a vertex in $\partial_o Y$. The set of external edges form the {\em Cheeger boundary} of $Y$ denoted by $\partial_{\ast}Y$. 

To any vertex $v$ in $\partial Y$ we can assign an external edge starting at $v$. This gives an injection from  $\partial Y$ to  $\partial_{\ast}Y$. On the other hand, there exist at most $2m$ edges in $\partial_{\ast}Y$ staring at $v$. This implies inequalities  $|\partial Y|\le|\partial_{\ast}Y|\le2m|\partial Y|$. The same arguments show that $|\partial_o Y|\le|\partial_{\ast}Y|\le2m|\partial_oY|$. 

By the {\em density} of a subgraph $Y$ we mean its average vertex degree. This concept was introduced in~\cite{Gu04}; see also~\cite{Gu21a}. It is denoted by $\delta(Y)$. A {\em Cheeger isoperimetric constant} of the subgraph $Y$ is the quotient $\iota_*(Y)=|\partial Y|/|Y|$. It follows directly from the definitions that $|\delta(Y)|+|\iota_*(Y)|=2m$. Indeed, each vertex $v$ has degree $2m$ in the Cayley graph $\Gamma$. This is the sum of the number of internal edges starting at $v$, which is $\deg_Y(v)$, and the number of external edges starting at $v$. Taking the sum over all $v\in Y$, we have $2m|Y|$, which is equal to $\sum\limits_v\deg_Y(v)+|\partial_{\ast}Y|$. Dividing by $|Y|$, we get the above equality.

By the density of the Cayley graph $\Gamma={\cal C}(G;A)$ we mean the number $\bdelta(\Gamma)=\sup\limits_Y\delta(Y)$, where $Y$ runs over all nonempty finite subgraphs. Analogously, the Cheeger isoperimetric constant of the group $G$ with respect to generating set $A$, is the number $\iota_*(\Gamma)=\inf\limits_Y\iota_*(Y)$. Clearly, $\bdelta(\Gamma)+\iota_*(\Gamma)=2m$.

In further sections, working with isoperimetric constants, we will write $\iota_*(G;A)$ instead of the above notation.

\subsection{Amenability and flows on graphs}
\label{amen}

Recall that the group is called {\em amenable} whenever there exists a finitely additive normalized invariant mean on $G$, that is, a mapping $\mu\colon{\cal P}(G)\to[0,1]$ such that $\mu(Z_1\cup Z_2)=\mu(Z_1)+\mu(Z_2)$ for any disjoint subsets $Z_1,Z_2\subseteq G$, $\mu(G)=1$, and $\mu(Zg)=\mu(gZ)=\mu(Z)$ for any $Z\subseteq G$, $g\in G$. One gets an equivalent
definition of amenability if only one-sided invariance of the mean is assumed, say, the condition $\mu(Zg)=\mu(Z)$ ($Z\subseteq G$, $g\in G)$. The proof can be found in \cite{GrL}.

We are not going to list all well-known properties of (non)amenable groups. It is sufficient to refer to one of modern surveys like~\cite{Sap14}. Just notice that all finite and Abelian groups are amenable. The class of amenable groups is closed under taking subgroups, homomorphic images, group extensions, and directed unions of groups. The groups in the closure of the classes of finite and Abelian groups under this list of operations are called {\em elemetary amenable} (EA). Also we need to say that free groups of rank $> 1$ are not amenable. 
\vspace{1ex}

We will often refer to the following F\o{}lner criterion~\cite{Fol}. Here we restrict ourselves to the case of finitely generated groups.

\begin{prop}
\label{fol}
A group $G$ with finite set of generators $A$ is amenable if and only if its Cheeger isoperimetric constant iz zero: $\iota_*(G;A)=0$.
\end{prop}

This holds for any finite set of generators. Equivalently, one can say that the density of the Cayley graph $\Gamma={\cal C}(G;A)$ (right or left) has its maximum value: $\bdelta(\Gamma)=2m$, where $m=|A|$.

In practice, to establish amenability of a group, it is sufficient to construct a collection of finite subgraphs $Y$ in the Cayley graph such that $\inf\limits_Y|\partial_{\ast}Y|/|Y|=0$. Such subsets of vertices are called {\em Folner sets}. Informally, this means that almost all vertices of these sets are internal.

If one expects that a group is not amenable, or the answer is unknown, it is reasonable to construct sets with the least possible  value of $\iota_*(G;A)$. We will discuss this approach later in Section~\ref{isop}.

In case when we are going to establish non-amenability of a group, it is useful to apply other known criteria.

Let $G$ be a group generated by a finite set $A$. A {\em doubling function} on $G$ is a mapping $\psi\colon G\to G$ such that

	a$)$ for all $g\in G$ the distance $\distan(g,\phi(g))$ is bounded from above by a constant $K > 0$,
	
	b$)$ any element $g\in G$ has at least two preimages under $\psi$.

\begin{prop}
\label{grom}
	A group $G$ with finite set of generators $A$ is non-amenable if and only if it admits a doubling functions.
\end{prop}

This criterion is often attributed to Gromov. An elegant proof of it based on the Hall -- Rado theorem can be
found in \cite{CGH}, see also \cite{DeSS}. Note that this property also does not depend on the choice of a finite generating set. 

An interesting partial case happens if the constant $K$ in the above definition equals $1$. In~\cite{Gu04} we proved that for a 2-generator group $G$, a doubling function with constant $K=1$ exists if and only if the density of the Cayley graph does not exceed $3$. Equivalently, one can say that the Cheeger isoperimetric constant is at least one: $\iota_*(G;A)\ge1$. 

We call a group $G$ {\em strongly non-amenable} (with respect to a given finite generating set $A$) whenever this inequality holds: $\iota_*(G;A)\ge1$. One of the goals of this paper is to find a nesessary and sufficient condition for this property for groups with any finite number of generators. This will be done in Section~\ref{ev}.
\vspace{1ex}

Another practical criterion for non-amenabilty can be stated in terns of flows on Cayley graphs. Ley us introduce the terminology for that.

A {\em flow} on a graph $\Gamma$ is a real-valued function $f\colon E\to\mathbb R$ such that $f(e^{-1})=-f(e)$. We say that $f(e)$ is the flow through the edge $e$. Given a vertex $v$, we define an {\em inflow} to it as a sum of flows through all edges with $v$ as a terminate vertex.

The following criterion in terms of flows is essentially known. It can be derived from the above criterion in terms of doubling functions.

\begin{prop}
\label{critfl}
Let $G$ be a group with finite generating set $A$, and let $\Gamma={\cal C}(G;A)$ be its Cayley graph. The group $G$  is non-amenable if and only if there exist constants $C > 0$ and $\varepsilon > 0$, and a flow $f$ on $\Gamma$ with the following properties:

a$)$ The absolute value of the flow through each edge is bounded: $|f(e)|\le C$ for all $e\in E$;

b$)$ The inflow into each vertex is at least $\varepsilon$.	
\end{prop}

Let us give a sketch of the proof. Suppose that the flow from the above statement exists. Take any finite nonempty subset $Y$ of $G$. The sum of inflows into all its vertices is at least $\varepsilon|Y|$. On the other hand, this inflow comes into $Y$ from the outside through the Cheeger boundary. Each of its edges can bring the flow at most $C$. This implies inequality $\varepsilon|Y|\le C|\partial_{\ast}Y|$. So $\iota_*(G;A)\ge\varepsilon/C > 0$. By Folner criterion, $G$ is non-amenable.

The converse can be established as follows. Let $G$ be non-amenable. By Proposition~\ref{grom}, let $\psi$ be a doubling function on $G$. For any $g\in G$ we choose the shortest path from $g$ to $\psi(g)$ in the Cayley graph. The lengths of all these paths are uniformly bounded by the constant $K$. To each of the paths we assign a unit flow through each of its edges. The sum of these flows gives the flow on the whole Cayley graph. The resulting flow through every edge is finite. Moreover, it is uniformly bounded. Namely, it does dot exceed the size of the ball of radius $K$ in the Cayley graph. This will be the constant $C$ for the flow, and $\varepsilon=1$ since each vertex has onle one image and at least two preimages. So the inflow to it will be at least $1$.

\subsection{Rooted binary trees and forests}

We add this short subsection to introduce some notation used in the paper. More details on the subject can be found in~\cite{Gu21a}.

Formally, {\em a rooted binary tree} can be defined by induction.

1) A dot $\cdot$ is a rooted binary tree.

2) If $T_1$, $T_2$ are rooted binary trees, then $(T_1\hat{\ \ }T_2)$ is a rooted binary tree.

3) All rooted binary trees are constructed by the above rules.
\vspace{1ex}

Instead of formal expressions, we will use their geometric realizations. A dot will be regarded as a point. It coincides with the root of that tree. If $T=(T_1\hat{\ \ }T_2)$, then we draw a {\em caret\/} for $\hat{}$ as a union of two closed intervals $AB$ 
(goes left down) and $AC$ (goes right down). The point $A$ is the {\em root} of $T$. After that, we draw trees for $T_1$, $T_2$ and
attach their roots to $B$, $C$ respectively in such a way that they have no intersection. It is standard that
for any $n\ge0$, the number of rooted binary trees with $n$ carets is equal to the $n$th Catalan number $c_n=\frac{(2n)!}{n!(n+1)!}$.

Each rooted binary tree has {\em leaves\/}. Formally, they are defined as follows: for the one-vertex tree
(which is called {\em trivial\/}), the only leaf coincides with the root. In case $T=(T_1\hat{\ \ }T_2)$, the set of leaves
equals the union of the sets of leaves for $T_1$ and $T_2$. In this case the leaves are exactly vertices of degree
$1$.

We will also need the concept of a {\em height\/} of a rooted binary tree. For the trivial tree, its height equals
$0$. For $T=(T_1\hat{\ \ }T_2)$, its height is $\hgt T=\max(\hgt T_1,\hgt T_2)+1$.

Now we define a {\em rooted binary forest\/} as a finite sequence of rooted binary trees $T_1$, ... , $T_m$,
where $m\ge1$. The leaves of it are the leaves of the trees. It is standard from combinatorics that the number
of rooted binary forests with $n$ leaves also equals $c_n$. The trees are enumerated from left to right and they
are drawn in the same way.

A {\em marked\/} (rooted binary) forest is a (rooted binary) forest where one of the trees is marked.

\section{Evacuation schemes}
\label{ev}

Let us introduce a new concept. We have an infinite group $G$ generated by a finite set $A$. Let $\Gamma={\cal C}(G;A)$ be its Cayley graph (right or left). To each vertex $v$ we assign an infinite path $p_v$  starting at $v$ in the Cayley graph. Suppose that there exists a constant $C$ such that each directed egde $e$ can participate in these paths at most $C$ times. In this case we say that the family $(p_v)_{v\in G}$ is an {\em evacuation scheme} on the Cayley graph $\Gamma$.

To be more precise: we claim that the total number of occurrences of each edge $e$ in paths of the form $p_v$ $(v\in G)$ does not exceed $C$.

In case if $C=1$ we say that we have a {\em pure evacuation scheme}. This will be the most interesting situation for us. We are going to give a criterion for existing such a scheme.

We start from the following easy observation.

\begin{prop}
\label{critev}
A group $G$ with finite generating set $A$ is non-amenable if and only if there exists an evacuation scheme on its Cayley graph.
\end{prop}

\prf Suppose that an evacuation scheme exists. Take any finite nonempty subset $Y$ of $G$. For each $v\in Y$, the infinite path $p_v$ must leave the set $Y$ through some edge. (Otherwise some edge of $Y$ will be involved in paths infinitely many times.) This edge belongs to the Cheeger boundary of $Y$. Each of these edges may evacuate at most $C$ vertices from $Y$ via their paths. Therefore, the total number of vertices that can be evacuated does not exceed $C|\partial_{\ast}Y|$. On the other hand, all vertices in $Y$ should go outside $Y$ at some moment. So we have inequality $|Y|\le C|\partial_{\ast}Y|$. Therefore, the Cheeger isoperimetric constant $\iota_*(G;A)=\sup\limits_Y|\partial_{\ast}Y|/|Y|\ge\frac1C > 0$. By Folner criterion, $G$ is non-amenable.

Let us also observe that in case of pure evacuation scheme one has $\iota_*(G;A)\ge1$ so the group is strongly non-amenable (with respect to $A$).
\vspace{1ex}

Let us prove the converse implication. If $G$ is non-amenable, we have a doubling function $\psi$ from Proposition~\ref{grom}. Let us choose the shortest path $s(g)$ in the Cayley graph from $g$ to $\psi(g)$. Notice that each egde $e$ occurs finitely many times in the paths of the form $s(g)$. Indeed, the lengths of these paths are uniformly bounded by the constant $K$. So if an edge $e$ occurs in a path $s(g)$, then the distance between $g$ and the initial vertex $v=\iota(e)$ of $e$ does not exceed $K$. Thus $g$ belongs to the ball of radius $K$ around $v$ in $\Gamma$. Its size can be roughly estimated as $C=1+2m+2m(2m-1)+\cdots+2m(2m-1)^{K-1}$. Therefore, $e$ occurs in at most $C$ paths of the form $s_g$. These paths are geodesic, so $e$ can occur only once in each of them. The total number of all these occurrences of $e$ will not exceed $C$.

Now let us construct an evacuation scheme. The group $G$ is countable. Let us enumerate all its elements: $g_1$, $g_2$, ... , $g_n$, \dots\,. By induction on $n$, we define the evacuation paths $p_v$ for each $v=g_n$. We regard $\psi$ as a binary relation on $G$, that is, a set of ordered pairs. At each step we will exclude some of them from $\psi$. Let $\psi_0=\psi$, and let $\psi_n$ be the set of pairs obtained after $n$ steps of the process. 

Given a binary relation (or a partial function) on the set of vertices, we define an {\em index} of a vertex with respect to it as the difference between the number of its preimages and the number of its images (we assume that both numbers are finite). We know that all indexes of vertices are positive with respect to $\psi$. 

We prove by induction on $n\ge1$ that the indexes of $v_1$, \dots, $v_{n-1}$ are nonnegative and the indexes of all other vertices are positive, with respect to $\psi_{n-1}$. This is true for $n=1$. At the $n$-th step, we take the vertex $v_n=w_0$. Its index is positive, so it has a preimage $w_1$ under $\psi_{n-1}$. Since the indexes are nonnegative, each vertex that has an image must have a preimage. So let $w_2$ be the preimage of $w_1$, and so on. After defining this infinite sequence of the $w_i$s, we delete all the ordered pais of the form $(w_{i+1},w_i)$ ($i\ge0)$ from $\psi_{n-1}$. The result is the relation $\psi_n$. Clearly, the index of $v_n$ decreases by $1$ and becomes nonnegative; the indexes of other vertices are without changes. So $\psi_n$ satisfies the desired condition.

Now we define an infinite evacuation path $p_v$ for $v=v_n$ as $s(w_1)^{-1}s(w_2)^{-1}\ldots$\,. This completes the inductive step.

After all the evacuation paths have been contructed, we observe that the number of occurrences of any edge $e$ in them does not exceed $C$ by the property of paths $s(g)$. This completes the proof.
\vspace{1ex}

Suppose that we have an evacuation scheme with constant $K$ on the Cayley graph $\Gamma={\cal C}(G;A)$. Let $Y$ be a finite nonempty subset of $G$. We know that each path $p_v$ ($v\in Y$) must leave $Y$ at some step. Hence there exists an initial segment $\bar{p}_v$ of $p_v$ that has its terminal point on $\partial Y$. This leads to the concept of an {\em evacuation scheme with constant $K$ on a finite subgraph}. This is a collection of paths in $Y$ of the form $\bar{p}_v$. This finite path starts at $v$ and ends on the inner boundary $\partial Y$. For each edge $e$ we claim that the total number of its occurrences in the paths $\bar{p}_v$ ($v\in Y$) does not exceed $K$.

\begin{rk}
\label{simp}
Having an evacuation scheme on $Y$, we can assume that paths $\bar{p}_v$ are simple (otherwise we can remove some loops). Also we can claim that if $e$ occurs in the evacuation paths, then $e^{-1}$ does not occur. Indeed, if $\bar{p}_v=p_1ep_2$, $\bar{p}_u=p_3e^{-1}p_4$, then one can replace these evacuation paths by $p_1p_4$, $p_3p_2$, respectively. 
\end{rk}

The following fact is an easy consequence of K\"onig's Lemma~\cite{Kon}.

\begin{prop}
\label{kon}
Let $G$ be an infinite group generated by a finite set $A$. An evacuation scheme on the Cayley graph $\Gamma={\cal C}(G;A)$ exists if and only if there exist evacuation schemes on all its finite nonempty subgraphs. The constants for both cases are the same.
\end{prop}

\prf The ``only if'' part is obvious. Let there exist evacuation schemes on all finite nonempty subgraphs of $\Gamma={\cal C}(G;A)$. (The constant is always assumed to be $K$.) Let us take a collection of balls $B(n)$ around the identity ($n\ge0$). For each of these balls we take all evacuation schemes on it. Ordered pairs of the form $\la B(n), {\cal E}\ra$ will form the set of vertices of an auxiliary directed graph, where ${\cal E}$ is an evacuation scheme on $B(n)$.

We put a directed edge from $\la B(n), {\cal E}\ra$ to $\la B(n+1), {\cal E}'\ra$ whenever for all $v\in B(n)$, the evacuation path of $v$ in $B(n)$ if the initial segment of the evacuation path of $v$ in $B(n+1)$. Clearly, the graph is locally finite since the group is finitely generated. It is also infinite, so we have an infinite chain in it according to K\"onig's Lemma. This infinite chain defines an infinite path starting at $v$ in $\Gamma$ (as an ascending ``union'' of evacuation paths for $v$ in the balls). All edges have a uniformly bounded number of occurences in the evacuation paths. So we have an evacuation scheme on $\Gamma$ with the same constant. Additionally, we can now assume that these evacuation paths are simple, and no edge can be involved in these paths together with its inverse.

The proof is complete.
\vspace{1ex}

Now we are ready to present the main result of this Section. 

\begin{thm}
\label{main2}
Let $G$ be an infinite group generated by a finite set $A$. A pure evacuation scheme on the Cayley graph $\Gamma={\cal C}(G;A)$ exists if and only if the group $G$ is strongly non-amenable with respect to $A$, that is, its Cheeger isoperimetric constant is at least one: $\iota_*(G;A)\ge1$.	
\end{thm}

\prf The ``only if'' part is easy. It has been already established in the proof of Proposition~\ref{critev}. Let us prove the converse. This can be done using Ford -- Fulkerson theorem~\cite{FF} in terms of maximal network flows and minimal cuts. However, we prefer a direct proof.

Suppose that $\iota_*(G;A)\ge1$. If there is no pure evacuation scheme on $\Gamma$, then there exists a finite nonenpty subset $Y$ without pure evacuation schemes on it, due to Proposition~\ref{kon}. Let us consider a partial pure evacuation scheme on $Y$ in the following sense. To some of vertices $v\in Y$ we assign simple paths in $Y$ from $v$ to a vertex in $\partial Y$. Each directed edge $e$ occurs in at most one of these paths. We choose this scheme in such a way that the subset $Y'$ of vertices $v$ for which we had assigned $p_v$, has maximal possible cardinality. Vertices from $Y\setminus Y'$ will be called {\em incomplete}. The set of them is nonempty. Further, we can assume that if $e$ occurs in paths $p_v$ ($v\in Y'$), then $e^{-1}$ does nor occur, according to Remark~\ref{simp}.

If a directed edge $e$ occurs in a path of the form $p_v$ ($v\in Y'$) then we call it {\em full}. A vertex $u$ is called {\em accessible} if there exists a path in $Y$ from some incomplete vertex $v$ to $u$ such that no directed edge of it is full. By $Z$ we denote the set of all accessible vertices. It contains all incomplete vertices so it is nonemty. Notice that vertices in $Z$ cannot belong to the inner boundary of $Y$. Otherwise we can add a path from an incomplete vertex to a vertex in the boundary of $Y$, extending the set of paths of the form $p_v$ and thus increasing the size of $Y'$.

Let us consider the Cheeger boundary of $Z$. All its egdes are inside $Y$ since $Z$ does not meet $\partial Y$. It follows that all directed edges in $\partial_{\ast}Z$ are full. Indeed, if an edge of this form is not full, then its terminal vertex is accessible since its initial vertex belongs to $Z$. However, the terminal vertex does not belong to $Z$ by definition of the Cheeger boundary. Thus we will have a contradiction.

Since each edge of $\partial_{\ast}Z$ is full, it belongs to a path of the form $p_v$. We claim that $v$ belongs to $Z$. Otherwise there is an edge in $p_v$ that goes from a vertex outside $Z$ to a vertex in $Z$. This egde will be inverse to an edge from $\partial_{\ast} Z$. This implies that mutually inverse edges are full contrary to our additional assumption.

The previous paragraph also shows that if a path of the form $p_v$ goes along an edge from the Cheeger boundary of $Z$, then it never returns to $Z$. Therefore, each edge from $\partial_{\ast}Z$ occurs in its own path of the form $p_v$, where $v\in Y'$. This gives an injection from  $\partial_{\ast}Z$ to $Z$. In particular, $|\partial_{\ast}Z|\le|Z|$. Besides, we know that all incomplete vertices belong to $Z$. At least one of them exists, and it is not an initial point of any path of the form $p_v$. Hence it is not in the image of the above injection. This implies strict inequality $|\partial_{\ast}Z| < |Z|$. Therefore, $\iota_{\ast}(G;A)\le|\partial_{\ast}Z|/|Z| < 1$ --- a contradiction.

The proof is complete.

\begin{rk}
\label{dou2}
Let us describe a relationship between pure evacuation schemes and generalizations of doubling maps. This gives a new equivalent characterization of the new concept.

Suppose that we have a pure evacuation scheme on the Cayley graph $\Gamma={\cal C}(G;A)$. For each edge $e$ that occurs in a path of the form $p_v$, we put an arrow on it in the opposite direction. This defines a multi-valued partial function $\psi$ from $G$ to itself. More formally, we take all ordered pairs of the form $\la\tau(e),\iota(e)\ra$, where $e$ runs over the above set of edges. The collection of these pairs is the binary relation $\psi$.

It follows directly from the construction that each vertex has a finite number of images and preimages under $\psi$, and the index of each vertex (the difference between the number of preimages and images) is always positive.

On the other direction, if we have a multi-valued partial function $\psi$ satisfying the above conditions, then we can extract from it a pure evacuation scheme on $\Gamma$. This is done exactly in the same way as in the proof of Proposition~\ref{critev}.
\end{rk}

\section{Isoperimetric constants for Thompson's group $F$}
\label{isop}

We define R.\,Thompson's group $F$ by the following infinite group presentation

\be{xinf}
\la x_0,x_1,x_2,\ldots\mid x_j{x_i}=x_ix_{j+1}\ (i<j)\,\ra.
\ee
This group was found by Richard J. Thompson in the 60s. In many papers the same group is defined via piecewise-linear functions. We do not use them in this paper so we prefer to be based on the combinatorial definition. We refer to the survey \cite{CFP} for details. See also~\cite{BS,Bro,BG}. We also refer to~\cite{GbS} where it was proved that $F$ is a diagram group; further development of this theory can be found in~\cite{GuSa99}. In~\cite{GuSa} we constructed a new normal form for elements of $F$ different from the standard one.

It is easy to see that for any $n\ge2$, one has $x_n=x_0^{-(n-1)}x_1x_0^{n-1}$ so
the group is generated by $x_0$, $x_1$. It can be given by the following presentation with two defining relations

\be{x0-1}
\la x_0,x_1\mid x_1^{x_0^2}=x_1^{x_0x_1},x_1^{x_0^3}=x_1^{x_0^2x_1}\ra,
\ee
where $a^b=b^{-1}ab$ by definition. Also we define a commutator $[a,b]=a^{-1}a^b=a^{-1}b^{-1}ab$
and notation $a\leftrightarrow b$ whenever $a$ commutes with $b$, that is, $ab=ba$.

An equivalent definition of $F$ can be given in the following way. Let us consider all strictly increasing continuous piecewise-linear functions from the closed unit interval onto itself. Take only those of them that
are differentiable except at finitely many dyadic rational numbers and such that all slopes (derivatives) are integer powers of $2$. These functions form a group under composition. It is is isomorphic to $F$. Another useful representation of $F$ by piecewise-linear functions can be obtained if we replace $[0,1]$ by $[0,\infty)$ in the previous definition
and impose the restriction that near infinity all functions have the form $t\mapsto t+c$, where $c$ is an integer.

The group $F$ has no free subgroups of rank $>1$. It is known that $F$ is not elementary amenable (EA). However, the famous problem about amenability of $F$ is still open. The question whether $F$ is amenable was asked by Ross Geoghegan in 1979; see~\cite{Ger87}. There were many attempts of various authors to solve this problem in both directions. We will not review a detailed history of the problem; this information can be found in many references. However, to emphasize the difficulty of the question, we mention the paper~\cite{Moore13}, where it was shown that if $F$ is amenable, then Folner sets for it have a very fast growth. Besides, we would like to refer to the paper~\cite{BB05} where the authors obtained an estimate of the isoperimetric constant of the group $F$ in its standard set of generators $\{x_0,x_1\}$. This estimate has not been improved so far.

Notice that if $F$ is amenable, then it is an example of a finitely presented amenable group, which is not EA. If it is not amenable, then this gives an example of a finitely presented group, which is not amenable and has no free subgroups of rank $>1$.  Note that the first example of a non-amenable group without free non-abelian subgroups has been constructed by Ol'shanskii \cite{Olsh}. (The question about such groups was formulated in \cite{Day}, it is also often attributed to von Neumann \cite{vNeu}.) Adian \cite{Ad83} proved that free Burnside groups with $m>1$ generators of odd exponent $n\ge665$ are not amenable. The first example of a finitely presented non-amenable group without free non-abelian subgroups has been constructed by Ol'shanskii and Sapir \cite{OlSa}. Grigorchuk \cite{Gri} constructed the first example of a finitely presented amenable group not in EA.
\vspace{2ex}

It is not hard to see that $F$ has an automorphism given by $x_0\mapsto x_0^{-1}$, $x_1\mapsto x_1x_0^{-1}$. To check that, one needs to show that both defining relators of $F$ in (\ref{x0-1}) map to the identity. This is an easy calculation using
defining relations of the group. After that, we have an endomorphism of $F$. Aplying it once more, we have the identity map. So this is an automorphism of order $2$.

Notice that $F$ has no non-Abelian homomorphic images~\cite{Bro}. So in order to check that an endomorphism of $F$ is a monomorphism, it suffices to show that the image of the commutator $[x_0,x_1]=x_0^{-1}x_1^{-1}x_0x_1=x_2^{-1}x_1=x_1x_3^{-1}$
is nontrivial.

Later we will work with the symmetric set $S=\{x_1,\bar{x}_1\}$ of generators, where $\bar{x}_1=x_1x_0^{-1}$. 
Obviously, it also generates $F$. Using Tietze transormations, one can show that $F$ has finite presentaion

\be{ab}
\la\alpha,\beta\mid\alpha^{\beta}\leftrightarrow\beta^{\alpha},\alpha^{\beta}\leftrightarrow\beta^{\alpha^2}\ra.
\ee
where $\alpha=x_1^{-1}$, $\beta=\bar{x}_1^{-1}=x_0x_1^{-1}$ (so $x_0=\beta\alpha^{-1}$).

Of course, from the symmetry reasons we know that $\beta^{\alpha}\leftrightarrow\alpha^{\beta^2}$
also holds in $F$. Therefore, it is a consequence of the two relations of (\ref{ab}). Moreover, one can check that
for any positive integers $m$, $n$ it holds $\alpha^{\beta^m}\leftrightarrow\beta^{\alpha^n}$ as a consequence
of the defining relations. More details on this presentation can be found in~\cite{Gu21a}.
\vspace{1ex}

Now let us discuss the question about isoperimetric constants for $F$ in variuos generating sets. First of all, for the standard generating set $A=\{x_0,x_1\}$, it was shown in~\cite{Gu04} that $\iota_*(F;A)\le1$. In the Addendum to
the same paper, there was a modification of the above construction showing that there are subgraphs with density
strictly greater than $3$, that is, $\iota_*(F;A) < 1$. A much stronger result was obtained in~\cite{BB05}. The authors constructed a family of finite subgraphs with density approaching $3.5$. That is, $\iota_*(F;A)\le\frac12$. This is the best known estimate for the isoperimetric constant for the present moment. There were many attempts of various authors to improve it, but they were not successful. So it is natural to formulate the following 

\begin{conj}
\label{cnj1}
The Cheeger isoperimetric constant for the group $F$ in its standard set of generators equals $\frac12$:
$$
\iota_*(F;\{x_0,x_1\})=\frac12.
$$
\end{conj}

Of course, if this conjecture is true, then $F$ is non-amenable. If it is not true, then disproving it would be a serious breakthrough towards amenability.
\vspace{1ex}

Now we will need to give a brief description of the Brown -- Belk construction. The reader can find more information in~\cite{Be04,BB05}. Notice that in~\cite{Gu21a} we gave quite a precise description of Brown -- Belk sets. Now we need to recall  the most important properties of them. The only thing we decided to change, is the following. Now we act by generators on the right, working with right Catley graphs. All marked forests we will deal with, now are {\em negative marked forests}, that is, mirror images of usual marked forests under a horizontal axis. So the roots of trees will be on the bottom, the leaves on the top. We will not recall again about these changes.

The following elementary lemma was proved in~\cite{Gu21a}. We will refer to it as the {\em symmetric property}.

\begin{prop}
\label{invlet}
Let $G$ be a finitely generated group and let $\Gamma={\cal C}(G,A)$ be its Cayley graph. Let $Y$ be a finite nonempty subgraph of $\Gamma$. Then for any $a\in A^{\pm1}$, the number of edges in the Cheeger boundary $\partial_{\ast} Y$ labelled by $a$ is the same as the number of edges in $\partial_{\ast}Y$ labelled by $a^{-1}$.
\end{prop}

Let $n\ge1$, $k\ge0$ be integer parameters. By $BB(n,k)$ we denote the set of marked forests that have $n$ leaves, and
each tree has height at most $k$. The group $F$ has a right partial action on this set. Namely, $x_0$ acts by
shifting the marker left if this is possible. The action of $x_1$ is as follows. If the marked tree is trivial,
this cannot be applied. If the marked tree is $T=(T_1\hat{\ \ }T_2)$, then we remove its caret and mark the tree
$T_1$. It is easy to see that applying $\bar{x}_1=x_1x_0^{-1}$ means the same, replacing $T_1$ by $T_2$ for the marked tree.

Now one can see that $x_1$ and $\bar{x}_1$ are totally symmetric. They generate $F$ so one can regard them as
one more pair of natural generators instead of the standard ones.

The action of $x_1^{-1}$ and $\bar{x}_1^{-1}$ is defined analogously. Namely, if the marked tree of a forest is
rightmost, then $x_1^{-1}$ cannot be applied. Otherwise, if the marked tree $T$ has a tree $T''$ to the right of
it, then we add a caret to these trees and the tree $T\hat{\ \ }T''$ will be marked in the result. Notice that
if we are inside $BB(n,k)$, then both trees $T$, $T''$ must have height $< k$: otherwise $x_1^{-1}$ cannot be
applied. For the action of $\bar{x}_1^{-1}$, it cannot be applied if $T$ is leftmost. Otherwise the marked tree $T$ has
a tree $T'$ to the left of it. Here we add a caret to these trees and the tree $(T'\hat{\ \ }T)$ will be marked in the
result. As above, both trees $T'$, $T$ must have height $< k$ to be possible to stay inside $BB(n,k)$.
\vspace{1ex}

These rules are important to understand the structure of the Cayley graph of $F$ (for many generating sets). They will be always assumed and used without references.
\vspace{1ex}

One can regard $BB(n,k)$ as a subset of vertices of the Cayley graph of $F$. This can be done for each of the three generating sets
$\{x_0,x_1\}$, $\{x_1,\bar{x}_1\}$, and $\{x_1,\bar{x}_1,x_0\}$.

For any fixed $k$, let $n\gg k$. Since any tree of height $k$ has at most $2^k$ leaves, any forest in $BB(n,k)$
contains at least $\frac{n}{2^k}$ trees. Therefore, if we randomly choose a marked forest, the probabililty for
this vertex of an automaton to accept both $x_0$, $x_0^{-1}$ approaches $1$. Now look at the probability to
accept $x_1	$. The contrary holds if and only if the marked tree is trivial. We may assume this tree is
not the rightmost one of the forest. Then we remove the trivial tree and move the marker to the right. As a
result, we obtain an element of $BB(n-1,k)$. The inverse operation is always possible. So the probability we
are interested in, equals $|BB(n-1,k)|/|BB(n,k)|$. It approaches some number $\xi_k$ as $n\to\infty$. If $k\to\infty$, then $\xi_k\to\frac14$. 

For the inverse letter $x_1^{-1}$, the probability not to accept it is the same by the symmetric property. We see that the number of  edges in the Cheeger boundary of $BB(n,k)$ approaches one half of the cardinality of this set. This means that the density of the set $BB(n,k)$ approaches $3.5$ and it gives the upper bound $\frac12$ for the Cheeger isoperimetric constant $\iota_{\ast}(F;\{x_0,x_1\})$.
\vspace{1ex}

We see that $Y=BB(n,k)$ accepts $x_1$ if and only if the marked tree is nontrivial. But the same holds for accepting $\bar{x}_1$. Both events hold with probability approaching $\frac14$. By symmetric property, the probabilities are the same for inverse letter. This looks like the Cheeger isoperimetric constant for the symmetric generating set $\{x_1,\bar{x}_1\}$ approaches $1$. However, this is not the case. In~\cite[Theorem 2]{Gu21a} we proved that the density of the Cayley graph in these generators strictly exceeds $3$. Equivalently, the Cheeger isoperimetric constant is strictly less than $1$. Now we are going to strengthen this result showing that even for the extended generating set $\{x_0,x_1,\bar{x}_1\}$, we get the isoperimetric constant strictly less that $1$.

\begin{thm}
\label{isopsym0}
The Cheeger isoperimetric constant of the Cayley graph of Thompson's group $F$ in extended symmetric generating set $\{x_0,x_1,\bar{x}_1\}$ is strictly less than $1$.
\end{thm}

According to Theorem~\ref{main2}, this means that there is no pure evacuation scheme on the Cayley graph of $F$ in these generators.
\vspace{1ex}

\prf We start from the set $Y=BB(n,k)$ as a subgraph in $\Gamma={\cal C}(F;A)$, where $A=\{x_0,x_1,\bar{x}_1\}$. 
Let us consider a marked forest of the form $\dots, T', T, T'', \dots$ where $T$ is marked. Suppose that $T$ is trivial.
Then the corresponding vertex of $Y$ does not accept $x_1$ as well as $\bar{x}_1$ (we cannot remove the caret). Additionally suppose
that both trees $T'$, $T''$ have height $k$. This means that we cannot apply neither $\bar{x}_1^{-1}$ (adding a
caret to $T'$ and $T$), nor $x_1^{-1}$ (adding a caret to $T$ and $T''$). So in this case we get an isolated
vertex in the Cayley graph with respect to $\{x_1,\bar{x}_1\}$.

Let us denote the set of the above vertices by $Y_0$. Clearly, for the extending generating set they always accept $x_0$ and $x_0^{-1}$. So they have degree $2$ in $\Gamma$. The key point of the proof of~\cite[Theorem 2]{Gu21a} was that the probablility $p$ to belong to $Y_0$ for a randomly choosen vertex in $Y$ is greater than some positive constant. Namely, we showed that $p=|Y_0|/|Y| > p_0=\frac1{260}$. 

Using the basic properties of the sets $BB(n,k)$, we find $k\gg1$ and $n=n(k)\gg1$ in such a way that the density of the corresponding subgraph will be arbitrarily close to $5$, namely, will exceed $5-\varepsilon$ for some small $\varepsilon > 0$. This is possible because the probability for a vertex in $BB(n,k)$ not to accept $x_1$ is $\frac14+o(1)$. The same for accepting $x_1^{-1}$, $\bar{x}_1$, $\bar{x}_1^{-1}$. The probaility not to accept $x_0^{\pm1}$ is $o(1)$.

Now we remove the subset $Y_0$ from $Y$ together with directed edges labelled by $x_0^{\pm1}$ incident to them. If $V$ is the set of vertices of $Y$ and $E$ the set of its directed edges, then we know that the density $\delta=\frac{|E|}{|V|} > 5-\varepsilon$. We removed $p|V|$ vertices and $4p|V|$ directed edges (each geomeric edge carries two directed edges). 

Hence the density of the new subgraph will be
$$
\frac{|E|-4p|V|}{|V|-p|V|}=\frac{\frac{|E|}{|V|}-4p}{1-p}=\frac{\delta-4p}{1-p} > \frac{5-\varepsilon-4p}{1-p}=5+\frac{p-\varepsilon}{1-p} > 5+\frac{p_0-\varepsilon}{1-p_0}.
$$
The sum of the density of the subgraph and its Cheeger isoperimetric constant is $2m$, where $m$ is the number of group generators. So the Cheeger isoperimetric constant will satisfy the following inequality:
$$
\iota_{\ast}(F;\{x_0,x_1,\bar{x}_1\}) < 1-\frac{p_0-\varepsilon}{1-p_0}=1-\frac{p_0}{2(1-p_0)}=\frac{517}{518}.
$$
whenever $\varepsilon=p_0/2$. This completes the proof.
\vspace{1ex}

Remark that we need to take very large sets of marked forests to make this effect possible. At least, these sets cannot be put inside the memory of computers to do some calculations with them. 
\vspace{1ex}

Notice that Theorem~\ref{isopsym0} also strengthens~\cite[Theorem 3]{Gu21a} where we esimated the size of outer boundaries of finite subgraphs of the same Cayley graph. 
\vspace{1ex}

Now one can ask the following question. What are natural systems of generators of $F$ for which we have a chance to find an evacuation scheme on them? That is, the isoperimetric constant will be at least $1$, as a conjecture? A good candidate for this property is the set $\{x_0,x_1,x_2\}$. Easy analysis of Brown -- Belk example leads to

\begin{prop}
\label{x012}
The Cheeger isoperimetric constant of the Cayley graph of Thompson's group $F$ in generating set $\{x_0,x_1,x_2\}$ is at least $1$.
\end{prop}

Indeed, accepting $x_2$ for a vertex in $BB(n,k)$ means that the tree to the right of the marker is nontrivial. The probablity of that is $\frac14+o(1)$, as in the case of $x_1$. So for external edges labelled by $x_1$, $x_1^{-1}$, $x_2$, $x_2^{-1}$, one can reach the value asymptotically equivalent to the number of vertices.

As in the case of the set $\{x_0,x_1\}$, there are no visible reasons how to improve the estimate. Thus we offer

\begin{conj}
\label{cnj2}
The Cheeger isoperimetric constant for the group $F$ in the set of generators $\{x_0,x_1,x_2\}$ equals $1$:
$$
\iota_*(F;\{x_0,x_1,x_2\})=1.
$$
\end{conj}

Let $Y$ be a finite nonempty subgraph of the Cayley graph $\Gamma={\cal C}(F;A)$. By $\nu_Y(a)$ we denote the number of vertices of the automaton $Y$ that do not accept $a\in A^{\pm1}$. This is the number of edges in the Cheeger boundary $\partial_{\ast}Y$ labelled by $a$. By symmetric property (Proposition~\ref{invlet}), we have $\nu_Y(a)=\nu_Y(a^{-1})$ for all $a$.

If Conjecture~\ref{cnj1} is true, then for any automaton $Y$ one has $|\partial_{\ast}Y|\ge\frac12|Y|$, that is, $\nu_Y(x_1)+\nu_Y(x_0)+\nu_Y(x_1^{-1})+\nu_Y(x_0^{-1})\ge\frac12|Y|$. This is equivalent to $\nu_Y(x_1)+\nu_Y(x_0)\ge\frac14|Y|$. Now recall that $x_0\mapsto x_0^{-1}$, $x_1\mapsto\bar{x}_1=x_1x_0^{-1}$ is an automorphism of $F$. Replacing $Y$ be its automorphic image, and the same for the edge labels, we obtain that $\nu_Y(\bar{x}_1)+\nu_Y(x_0^{-1})\ge\frac14|Y|$. Hence $\nu_Y(x_1)+\nu_Y(\bar{x}_1)+2\nu_Y(x_0)\ge\frac12|Y|$, so
$$
\nu_Y(x_1)+\nu_Y(\bar{x}_1)+2\nu_Y(x_0)+\nu_Y(x_1^{-1})+\nu_Y(\bar{x}_1^{-1})+2\nu_Y(x_0^{-1})\ge|Y|.
$$
The left-hand side of the above inequality is the cardinality of the Cheeger boundary of $Y$ with respect to the generating multi-set $\{x_1,\bar{x}_1,x_0,x_0\}$, where $x_0$ is doubly repeated. So we have one more conjecture, implied by Conjecture~\ref{cnj1}.

\begin{conj}
\label{cnj3}
The Cheeger isoperimetric constant for the group $F$ in the multi-set of generators $\{x_1,\bar{x}_1,x_0,x_0\}$ is at least $1$:
$$
\iota_*(F;\{x_1,\bar{x}_1,x_0,x_0\})\ge1.
$$
That is, there exists a pure evacuation scheme on the Cayley graph of $F$ for this generating set.
\end{conj}

If we compare this with Theorem~\ref{isopsym0}, then we see that the difference is only doubling the generator $x_0$, which is equivalent to allowing us to travel by egdes labelled by $x_0$ twice in the evaquation scheme for $\{x_1,\bar{x}_1,x_0\}$. This looks quite surprising since we know from the Brown -- Belk construction that $x_0^{\pm1}$ is almost always accepted by $BB(n,k)$.

It turns out that our second Conjecture also implies the above fact.

\begin{thm}
\label{c1c2c3}
Each of the Conjectures~\ref{cnj1} and~\ref{cnj2} implies Conjecture~\ref{cnj3}.
\end{thm}

\prf We have already proved that Conjecture~\ref{cnj1} implies Conjecture~\ref{cnj3}. Let Conjecture~\ref{cnj2} be true. Then there exists a pure evacuation scheme on the Cayley graph of $F$ in generators $\{x_0,x_1,x_2\}$. Conjugation by $x_0^{-1}$ takes the generators to $\{x_0,x_0x_1x_0^{-1}=x_0\bar{x}_1,x_1\}$. Each evacuation path going through an egde labelled by $x_0$, $x_1$, $x_2$ respectively, will now go through edges labelled by  $\{x_0,x_0\bar{x}_1,x_1\}$; the same for inverse labels. So we obtain an evacuation scheme on the Cayley graph for the set $\{x_0,x_1,\bar{x}_1\}$, where edges with label $x_0^{\pm1}$ can occur twice. This is equivalent to the pure evacuation scheme on the Cayley graph for the multi-set $\{x_1,\bar{x}_1,x_0,x_0\}$.

The proof is complete.
\vspace{1ex}

So we have the following alternative: either we confirm Conjecture~\ref{cnj3} proving non-amenability of $F$, or we disprove both Conjectures~\ref{cnj1} and~\ref{cnj2}, what will be a breakthrough in the direction to prove amenability of the group.

\end{document}